\documentclass{aims}
\usepackage[titletoc,toc,page]{appendix}
\usepackage{amsmath}
\usepackage{ amssymb}
  \usepackage{paralist}
  \usepackage{graphics} %% add this and next lines if pictures should be in esp format
  \usepackage{epsfig} %For pictures: screened artwork should be set up with an 85 or 100 line screen
\usepackage{graphicx}  
\usepackage{epstopdf}%This is to transfer .eps figure to .pdf figure; please compile your paper using PDFLeTex or PDFTeXify.
 \usepackage[colorlinks=true]{hyperref}
\hypersetup{urlcolor=blue, citecolor=red}
\usepackage{pictexwd,dcpic}    

\usepackage{subfig} % extra packages for figures and Taylor diagrams
\usepackage{color}
\usepackage{caption}
\usepackage{cancel}

\usepackage{nomencl} % nomenclature package. Useful for listing the symbols used in the main text
\makenomenclature

\newcommand {\R} {\mathbb {R}}

\newcommand {\C} {\mathbb {C}}
\newcommand {\Z} {\mathbb {Z}}

\newcommand {\Pj} {\mathbb {P}}

\newcommand {\no} {\noindent}

\def\currentvolume{X}
 \def\currentissue{X}
  \def\currentyear{200X}
   \def\currentmonth{XX}
    \def\ppages{X--XX}
     \def\DOI{10.3934/xx.xx.xx.xx}

\newtheorem{theorem}{Theorem}[section]

\theoremstyle{definition}
\newtheorem{definition}[theorem]{Definition}

\newtheorem{example}{Example}

\title[Bridges] %Use the shortened version of the full title
      {Bridges between subriemannian geometry and algebraic geometry: now and then}

\author[Alex Castro, Wyatt Howard, Corey Shanbrom]{}

 % Fixed
  % include more

 \email{adecastr@imperial.ac.uk, wyahoward@gmail.com, corey.shanbrom@csus.edu}
\thanks{A.C. and W. H. thank CAPES the partial support from Science without Frontiers BEX 11784-13-0 and Jovens Talentos 062/2013 respectively}

\begin{document}

%The abstract of your paper
\begin{abstract}
We consider how the problem of determining normal forms for a specific class of nonholonomic systems leads to various interesting and concrete bridges between two apparently unrelated themes. Various ideas that traditionally pertain to the field of algebraic geometry emerge here organically in an attempt to elucidate the geometric structures underlying a large class of nonholonomic distributions known as Goursat constraints. Among our new results is a regularization theorem for curves stated and proved using tools exclusively from nonholonomic geometry, and a computation of topological invariants that answer a question on the global topology of our classifying space. Last but not least we present for the first time some experimental results connecting the discrete invariants of nonholonomic plane fields such as the RVT code and the Milnor number of complex plane algebraic curves. 
\end{abstract}

\maketitle

% Enter the first author's name and address:
\centerline{\scshape Alex L Castro}
\medskip
{\footnotesize
% please put the address of the first author
 \centerline{Department of Mathematics}
   \centerline{Imperial College, London}
   \centerline{180 Queen's Gate, London SW7 2AZ, UK}
} % Do not forget to end the {\footnotesize by the sign }

%i am commenting this out
\medskip

\centerline{\scshape Wyatt Howard$^1$ and Corey Shanbrom$^2$}
\medskip
{\footnotesize
 % please put the address of the second  and third author
\centerline{$^1$ Mathematics and Computer Science, Santa Clara University, Santa Clara, CA, USA $\&$ }
\centerline{ Department of Mathematics, De Anza College, Cupertino, CA, USA}  
\centerline{$^2$ Department of Mathematics and Statistics, California State University, Sacramento, CA, USA}
%\centerline{Springfield, MO 65810, USA}
}

\bigskip

% The name of the associate editor will be entered by an editorial staff
% "Communicated by the associate editor name" is not needed for special issue.
 \centerline{(Communicated by the associate editor name)}

%%
%% LaTeX can automatically make a table of contents.  This is done by
%% uncommenting the following:
%%

\tableofcontents

%%
%%  To enter text is easy.  Just type it.  A blank line starts a new
%%  paragraph. 
%%

%%%%%%%%%%%%%%%%%%%%%%%%%%%%%%%%%%%%%%

%\nomenclature{$S(k)$}{$k$-step Semple tower }%
%\nomenclature{$\Delta_k$}{Plane field canonically associated with $S_k$ }%
%\nomenclature{$H^*(\bullet)$}{Cohomology ring of a variety }
%\nomenclature{$\Phi ^{n}$}{A symmetry of the $n$-th step Semple Tower }
%\nomenclature{$\omega$}{The $RVT$ code of a point, i.e. a string of $R$'s, $V$'s,$T$'s, and $L$'s }
%\nomenclature{$\Gamma(p)$}{Curves integral to the distribution at the point $p$ }
%\nomenclature{$\mbox{sgv}(p)$}{Small growth vector of a the distribution $\Delta$ at $p$}
%\printnomenclature 

%\noindent {\bf Terminology.} Brief explanation of the RVT code here. 
%{\color{blue} Corey, Wyatt: what do you guys think?}
%{\color{cyan} Don't feel like it is super necessary.  Especially since it will take up more space and we summarize these concepts in the relevant sections.}

%\begin{table}[t!]
%\caption{Outline of Common Notation}
%\begin{tabular}{|c|c|c|}
%\hline
%   Notation & Description  \\
%\hline
%   $S(k)$ & $k$-step Semple Tower  \\ \hline
%   $\Delta_{k}$ & Plane field canonically associated with $S_{k}$   \\ \hline
%   $H^{\ast}(M)$ & Cohomology ring of the variety $M$  \\ \hline
%   $\Phi ^{n}$ & A symmetry of the $n$-th step Semple Tower \\ \hline
%   $\omega$ & The $RVT$ code of a point; string of $R$'s, $V$'s,$T$'s, and $L$'s 
%      \\ \hline
%\end{tabular}
%\label{tab:notation}
%\end{table}

%%%%%%%%%%%%%%%%%%%%%%%%%%%%%%%%%%%%%%

%The title of your section 1
\section{Introduction}

%\noindent {\color{blue} Write down the main message: Novelty of this paper?} 

%\medskip

%\noindent {\color{red} Wyatt's write up on Chern has been included. Now it remains to include bibliography and our other contributions.} 

\begin{figure}% engineering solution
 %\label{fig:cartrailer}
  \centering
  \def\svgwidth{.6 \columnwidth}
   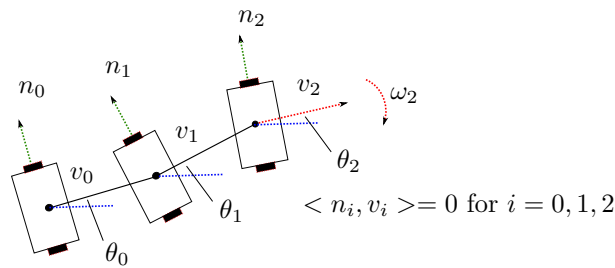
    \caption{The car with trailers attached. The rate of change in the steering direction is denoted by $\omega_2$.}
    \label{fig:cartrailer}
\end{figure}

One of the simplest idealized constraints one considers in nonholonomic mechanics is the {\em skate} or {\em no-slip} condition: % {\color{magenta} What are $x, y, \theta$?   }
\begin{equation}\label{eq:noslip}
-\sin(\theta) dx + \cos(\theta) dy = 0.
\end{equation}
where $(x,y)$ are cartesian coordinates in the plane and $\hat{v} = (\cos(\theta),\sin(\theta))$ is the steering direction of the car -- point along its axle. 
%{\color{magenta} Should 'normal' be 'parallel' here or am I crazy?  What does 'longitudinal' mean? What is $\omega_2$ in the figure? }
Stated plainly, the no-slip condition states that the wheels of the car are only allowed to roll along the road and, hence, no slipping occurs in the direction normal to the linear velocity of the car.  Another 
closely related constraint is that which an airport luggage cart is subject to (Figure \ref{fig:cartrailer}). In differential geometry such differential constraint is an example of a nonholonomic constraint, and it 
is well known that it does not possess integrable surfaces, i.e., it fails to satisfy the Frobenius condition of integrability (\cite{arnold1}, Appendix 3). 
%{\color{magenta} Is this reference necessary?  Are no-slip constraints really the exact same thing as contact distributions?  I think not.} 

%%%%%%%%%%%%%%%%%%%%%%%%%%%%%%
% Car with trailers picture

We now present to the reader a series of mathematical miracles related to the contact distribution above, and her close relatives the so-called {\em Goursat distributions}. 

%\medskip 

\begin{figure}
 \raggedright
    \subfloat[Prolongation of a distribution]{\label{fig:prolongdist} \def\svgwidth{.3 \columnwidth}
     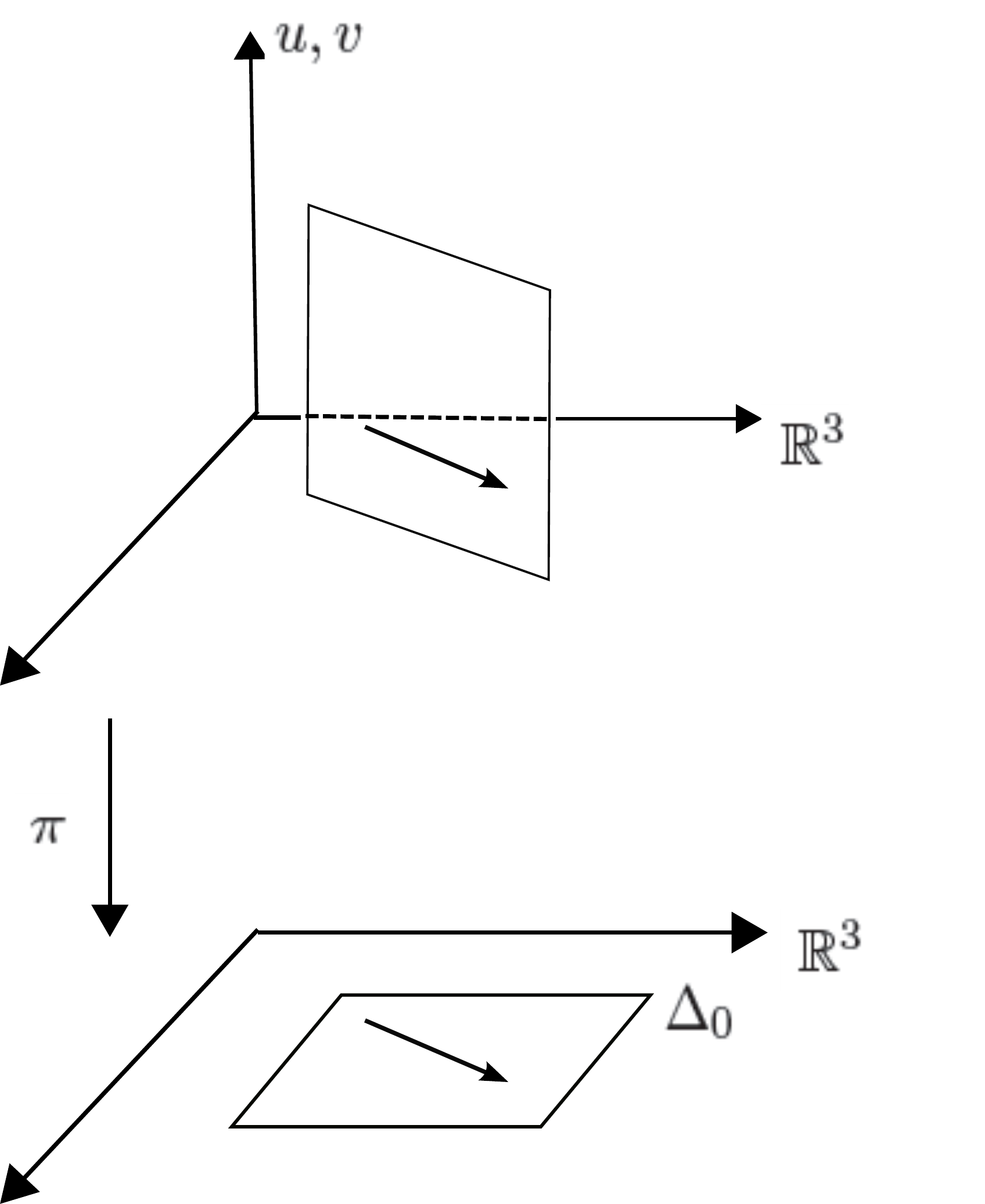}
    \subfloat[Prolongation of a diffeomorphism]{\label{fig:prolongdiff} \def\svgwidth{.3 \columnwidth}
     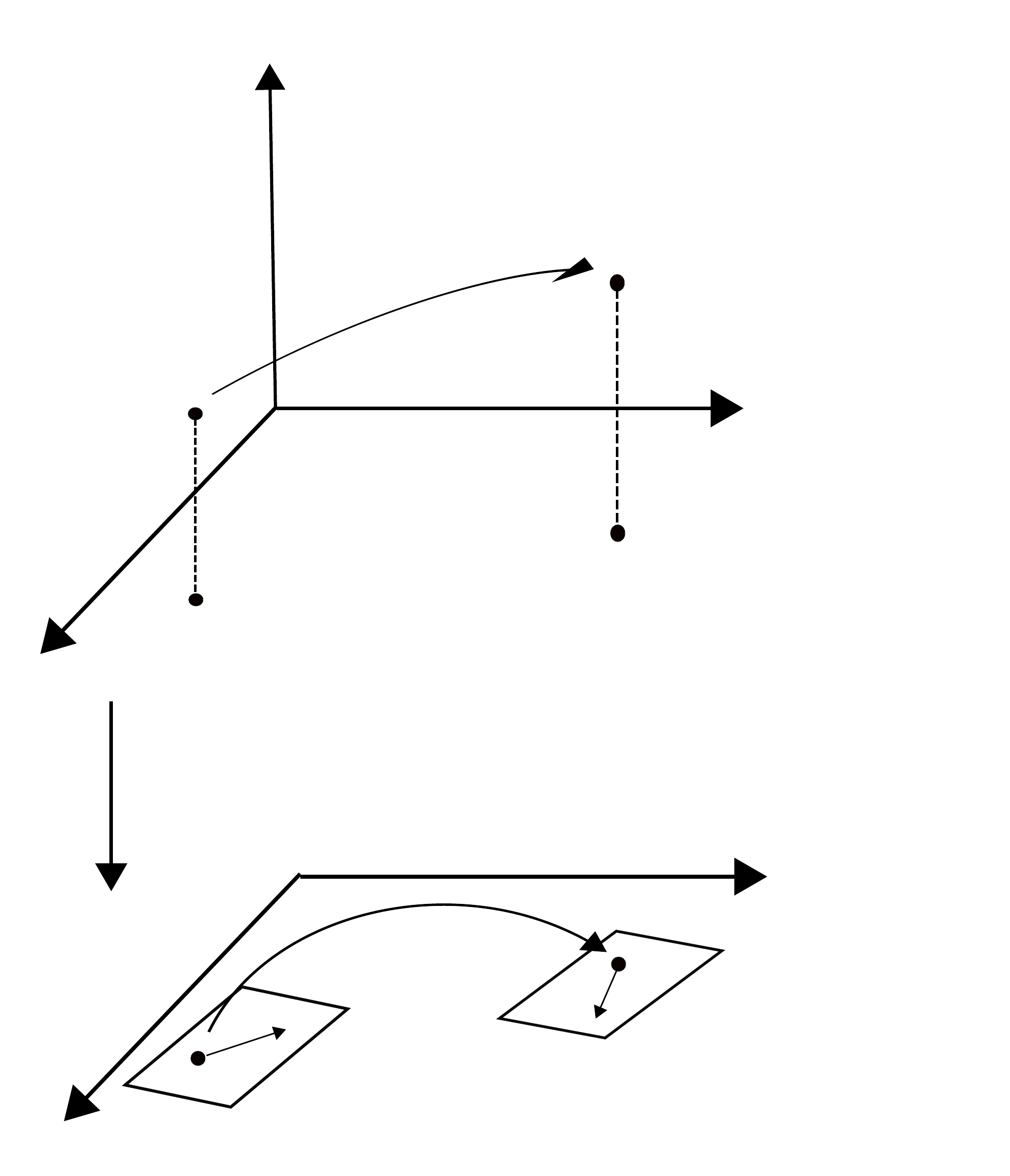}
     \subfloat[Prolongation of a curve]{\label{fig:prolongcur} \def\svgwidth{.3 \columnwidth}
     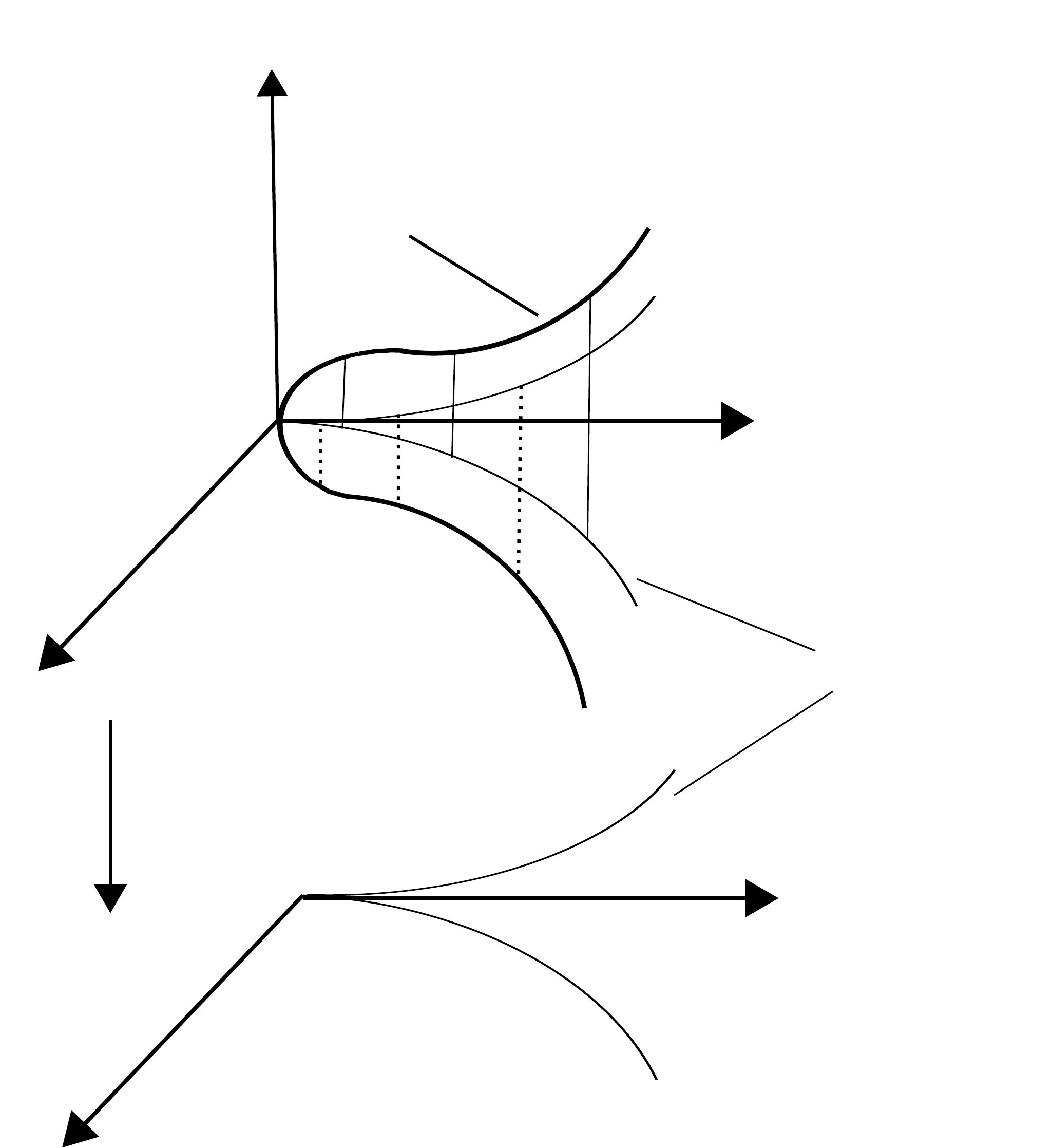}
     \caption{Prolongations}
     \label{fig:prolongation}
  \end{figure}

To begin, consider the projectivization of $T\mathbb{R}^2$ which we will denote by $S(1)$. We denote $\mathbb{R}^2$ by $S(0)$. It is a simple exercise to show that $S(1)$ is diffeomorphic to $
\mathbb{R}^2\times S^1$. The projectivization lifts various objects canonically defined in $\mathbb{R}^2$, including the tautological 2-plane field, morphisms and curves. In Figure \ref{fig:prolongation} we 
schematize the way various objects are lifted from $S(0)$ to $S(k)$.

%%%%%%%%%%%%%%%%%%%%%%%%%%%%%%%%%
% Figures for the Cartan prolongation

%%%%%%%%%%%%%%%%%%%%%%%%%%%%%%%%%

%{\color{magenta}{ I moved Fig 2 up a page, and removed two skips here}}

%\medskip

We can coordinatize $S(1)$ by using projective coordinates on $\mathbb{R}^2$: $(x,y,[dx:dy])$ where the third entry represents the projective coordinates of a vector in $\mathbb{R}^2$. Let $\pi:S(1) \to S(0)$ denote the canonical projection of the bundle just defined. Fix $x\in S(0)$ and consider a line $l \subset T_x S(0)$. Then the plane field 

$$\Delta_1 (x,[l]) \doteq D\pi^{-1}(l).$$ 

%{\color{magenta} The notation above is strange to me, and the first half of the first sentence below does not make sense to me.  Is $\theta=[dx:dy]$ ?}

The line $l$ has implicit representation in the plane given by $- \sin(\theta) dx + \cos(\theta) dy = 0$. In local coordinates $[dx:dy] \cong \tan(\theta)$ 
and the inverse image in $S(1)$ of this line under the tangent map $D\pi$ is the plane spanned by the vectors $\left\{ \frac{\partial}{\partial \theta}, \cos(\theta) \frac{\partial}{\partial x} + \sin(\theta) 
\frac{\partial}{\partial y}\right\}$. We have thus obtained the contact plane field (a.k.a skate constraint) from a geometric procedure known as {\bf Cartan prolongation}. This {\em operation} can be iterated, 
and at each step we projectivize the plane field obtained in the previous step and the output is a tower of fiber bundles known as the Semple Tower (\cite{kennedy1}):

$$\dots \to S(n) \to S(n-1) \to \dots \to S(2) \to S(1) \to S(0).$$ 

%\noindent {\color{magenta} From Corey: This is a good place to introduce the word 'Monster'.} \\
%{\color{blue} To Corey: any suggestions on how to introduce the word monster here?}

The Semple Tower has been rediscovered as the Monster Tower (\cite{montgomerygeo},\cite{montgomerypoints}).
One can assign to each point a word in the letters R,V, and T, known as an \emph{RVT code}, which is an invariant for the diffeomorphism group action.
Each consecutive level $S(k)$ is endowed with a plane field $\Delta_k$ which is a Goursat Flag and a geometric model for the configuration space of a airport luggage convoy (\cite{murray}, 
\cite{montgomerygeo}).  Some of these distributions had already been named within the literature of differential geometry; see {\cite{giarokumperaruiz}}.
%{\color{cyan} 'baptized' sounds odd. 'christened' sounds better or 'Some of these distributions were already known in the differential geometry literature by different names'.  Also, there is no Appendix in 
% this paper.}. 

One can generalize this construction to obtain a similar tower starting with a $d$-dimensional base manifold, but at each step one obtains a tower of fiber bundles with fibers isomorphic to $\mathbb{R}
P^{d-1}$ (or $S^{d-1}$ if orientation needs to be taken in into account). Each level of this tower will be also denoted by $S(k)$ without explicit mention to the base manifold. A word of caution here: 

$$\mbox{dim}(S(k)) = d + (d-1)k,$$ 

\noindent where $d$ is the dimension of the base manifold. Correspondingly, each $S(n)$ is equipped with a $d$-dimensional plane field. 

This tower can be thought of as the configuration space of a mechanical articulated arm \cite{pelletier}. There is also a mechanical resemblance between the articulated arm system just mentioned and the 
robotic snake model of Hausmann and Rodriguez \cite{hausmannrodriguez}, though we have to fix the position of the snake's tail. A simple geometric computation points out that for when the snake is completely stretched out the rank drops and it is not clear what the generic local behavior of this nonholonomic constraint is. Is this a Goursat distribution? Or a 
product constraint? By {\em product constraint}, we mean a plane field which contains factors as a product of a nonholonomic factor and trivial flat factor: $\Delta \times \mathbb{R}^k$. See section 5.3 of (\cite{montgomerygeo}). Hausmann and Rodriguez determined the reachable sets for 
certain generic configurations of the snake, but have not discussed in detail the nature of the nonholonomic 
distributions  (\cite{hausmannrodriguez}) . 

%{\color{cyan} Alex: We worked through the calculations of Hausmann and Rodriguez and showed that their distribution drops rank and can't in general give a Goursat distribution or the Monster Tower.  
%So, the Hausmann Rodriguez citation seems misleading in the beginning part of the paragraph.  The citation should be moved to the italicized part and state more clearly that the distribution arising from 
% the Snake Charmer Algorithm is not in general a Goursat distribution.}

The multi-flags distributions will be one of the protagonists of our note, and by a theorem by Y. Shibuya and K. Yamaguchi {(\cite{shibuyayamaguchi})} these generalized towers realize all 
of the so-called {\em Goursat multi-flags}. 
%{\color{magenta} Do we really need boldface above?}

\begin{quote}
Our take home message to the reader is that serendipity is abundant when it comes to Goursat flags and multi-flags.
\end{quote}

%%%%%%%%%%%%%%%%%%%%%%%%%%%%%%%%%

\section{Bridge 1: nonholonomic geometry and curve singularities}

A first mathematical miracle in the study of Goursat flags or multi-flags is their connection with the singularity theory of smooth of analytic maps \cite{arnold3}. The fundamental notion underlying this connection is the {\em Cartan prolongation} already alluded to in the introduction. From the point of view of normal form theory, there is an action of the pseudogroup of local diffeomorphisms at the base $S(0)$. Different orbits will correspond to different normal forms. Given a point $p\in S(k)$ one associates to it $\Gamma(p)$ which is the set of all smooth curves through $p$ that admit a nontrivial projection back to the base $S(0)$. It can be argued that $p\sim q$ ($\sim$ means equivalent under the group action) implies $\Gamma(p) \sim \Gamma(q)$. By symmetry, we fix the base point to be the origin and one can act on the fiber $\pi_k^{-1}(0)$, $\pi_k$ being the canonical projection $S(k)\mapsto S(0)$. On this set of curves we act with the diffeomorphism group and using standard techniques of singularity theory of curves, different curve orbits (normal forms) will give rise to different normal forms of plane fields. The first successful step taken in this direction was documented in \cite{montgomerypoints}. 

%{\color{blue} Corey: more understandable?}

% {\color{magenta}I don't understand the second to last sentence above.}

The task of labeling the orbits is sequential and it works by stages. In \cite{castro} we described in detail the general tools for this task  for the Semple Tower with base $S(0) = \mathbb{R}^3$, though the construction and tools work for general base manifolds. There are two main approaches for determining the orbits: 
\begin{enumerate}
\item the curve approach, and 
\item the isotropy method. 
\end{enumerate}
Both methods are rather elementary and perform well in lower dimensions, depending mostly on the combinatorial or projective geometry of the problem at hand. Limitations to both methods are mostly computational in nature, since the amount of bookkeeping grows exponentially as one goes up the tower. Either approach is suitable to computer algebraic systems as pointed out in (\cite{castro}). 

%The {\em curve method} consists in determining all classes of singular curves within each set $\Gamma(p)$. This task of determining the orbits on each $\Gamma(p)$ can be a Herculean task sometimes though. But curves can be described by essentially combinatorial invariants that simplify the analysis of their singularities. 

%One simple and illustrative example is the following. A certain prolongation of the contact distribution $du - p dx = 0$ will give rise to the following 2-plane field in $S(2)$:
%$$\left\{ \begin{array}{c} du - p dx = 0 \\ dx - p_2 dp = 0.\end{array}\right.$$ We say a curve is tangent to a {\bf critical direction} if its projection down to the base is constant. 
%The {\em critical directions} here are spanned by $\frac{\partial}{\partial p}, \frac{\partial}{\partial p_2}$. Take $p(t) = a t + O(2)$ and $p_2(t) = b t + O(2)$ with $ab\neq 0$. Then, 
%$$x(t) \propto t^2 + O(3),$$ and
%$$u(t) \propto t^3 + O(4).$$ This is a $A_2$ type of singularity. Any such curve is right-left equivalent to a cusp, i.e. $(t^2,t^3)$. Generically, a distribution only has curves which are prolongations of straight lines. Plane fields that fail to fulfill this condition are called {\bf singular}. 
%\medskip 

\subsection{Glimpses of algebraic geometry: Nash and Enriques} 

Many features from classic algebraic geometry are easily transported to the Semple Tower if one thinks in terms of Nash blow-ups instead of quadratic transformations. For a modern reference on the subject see \cite{lejeune}. 

Cartan prolongations coincide with Nash blow-ups in the analytic category as explained in \cite{castro5} and permitted us to make a case for the analogies between the problems in normal form theory for Goursat flags and some corresponding problems in enumerative geometry (\cite{kennedy1}). A technical but rather important result in the classification problem of Goursat multi-flags consists of the following: 

\begin{theorem}[\cite{castro5}, Appendix B]
Any well parametrized curve germ $c:I\to S(0)$ that is singular (i.e., $c'(0)=0$) becomes regular (i.e. smooth and touching only regular points) after a finite number of Cartan prolongations. 
\end{theorem}

%{\color{magenta}From Corey: Note mention of RVT below.  Also $c^{(k)}$ notation is unfamiliar, as is word 'attained'.  Consider deleting the discussion of the proof altogether.}\\ 
%{\color{blue} To Corey: What do you suggest?} 

This theorem is crucial in defining the RVT code of a point in $S(k)$ in terms of Cartan prolongations of curves in the base $S(0)$. Let us use the notation $c^{(k)}$ for the $k$-th iteration of the Nash blow-up of a curve $c(t)$. Consider the list of points $\{c^{(k)}(0)\}$ obtained by evaluating the consecutive prolongations at $t=0$. The proof consists of first showing that the set of points $c^{(k)}(0)\in S(k)$ contains only a finite number of critical points, and once we surpass the last critical point the curve become necessarily smooth. Otherwise it would have to be ill-parametrized. This is equivalent to the desingularization theorem stated in {\cite{lejeune}} but now formulated and proved using the language and tools of nonholonomic geometry. 

As a corollary of the regularization theorem we obtain a nonholonomic version of the classical Enriques theorem about multiplicities of consecutive prolongations of a well-parametrized curve, though originally it was worded in terms of quadratic transformations and proximity relations. By {\em multiplicity} of a singularity we mean the first non-zero jet (if we mod out constants due to a specific choice of chart). This definition is suitable to parametrized curves and is independent of the coordinate chart by di Bruno's formula.  %\footnote{\color{blue} http://mathworld.wolfram.com/FaadiBrunosFormula.html}. 
We will exchangeably use the term multiplicity for either points or curves. Using special charts known as extended Kumpera-Ruiz coordinates (see \cite{castro2}), the nested structure of Semple Towers (submanifolds of the original base manifold generate subtowers), and the regularization theorem above, we can prove the following: 

%{\color{magenta} Don't know what 'Semple tower theories' are above.}

\begin{definition}[\cite{castro5}, Appendix B]
Let $p \in S(k)$.  If a point $q$ {satisfies}
\begin{itemize}
\item  $q$ is in the fiber above $p$,  \emph{or}
\item  $q$ can be be reached by a prolongation of a vertical curve curve through $p$,
\end{itemize}
then we say that $p$ and $q$ are {\em adjacent points.}
Points which are connected this way will form a graph (in fact, a tree) with seed $p$. The adjacency condition will be denoted by $q \to p$.
\end{definition}
There is a simple relation between multiplies of adjacent points:
\begin{theorem}
One has
$$\mbox{mult}(p) = \sum_{q\to p} \mbox{mult}(q).$$
\end{theorem}  
This a classic result in enumerative geometry attributed to F. Enriques, and in our context it restricts the class of singular points that be reached from a given point via Cartan prolongation. The proof is again based purely on the local geometry of the nonholonomic fields in the Semple tower. 
%%%%%%%%%%%%%%%%%%%%%%%%%%%%%%%%%

\subsection{Puiseux numbers and growth vectors}

In \cite{Shanbrom}, we gave a formula for the Puiseux characteristic of an analytic plane curve germ which represents a Goursat distribution germ with prescribed small growth vector.  

Given a Goursat distribution $D$ on a manifold $M$, consider the sequence $D_{i}=[D, D_{i-1}]+D_{i-1}$, where $D_{0}=D$.  
Then there exists an $r$ such that $D_{r}=TM$.  
For each $p \in M$, we define the \emph{small growth vector at p} to be the integer valued vector
    \[sgv(p)= \big(\text{dim} D_{0}(p),\ \text{dim} D_{1}(p), \dots ,\ \text{dim} D_{r}(p)=n\big).  \]
   
    The \emph{derived vector} of a Goursat germ consists of the multiplicities of the entries in the small growth vector.
    For a Goursat distribution, the dimensions of the sequence $D_{i}$ grow by at most one at a time, so the multiplicities are nonzero and from the list of multiplicities we may recover the original small growth vector.

For a well-parametrized, non-immersed plane curve
\[ \gamma(t)=(t^m, \sum_{k \geq m} a_k t^k) \]
the Puiseux characteristic is defined as follows.
Let $\lambda_0=e_0=m$.  Then define inductively 
\[ \lambda_{j+1} = \min \{ k \ | \ a_k \neq 0,\ e_j \nmid k\}, \quad \quad e_{j+1}= \text{gcd}(e_j, \lambda_{j+1}) \]
until we first obtain a $g$ with $e_g=1$.  Then the vector $[\lambda_0 ; \lambda_1, \dots, \lambda_g]$ is called the \textit{Puiseux characteristic} of $\gamma$.
The Puiseux characteristic is the fundamental invariant in the singularity theory of plane curves.  In \cite{W}, Proposition 4.3.8 shows that it is equivalent to at least seven other classical invariants.
    
In short, \cite{Shanbrom} provided the dashed arrow in the following diagram:

%{\color{blue} To Corey: I couldn't generate the diagram you proposed. I may need your help here.} \\

\medskip 
    
    \[ \begindc{\commdiag}[60]
    \obj(0,1)[SGV]{\{SGV\}}
    \obj(1,1)[RVT]{\{RVT\}}
    \obj(1,0)[PC]{\{PC\}}
    %\mor{SGV}{RVT}{}[+1, 10] next two are replacing this one
    \mor{SGV}{RVT}{}
    \mor{RVT}{SGV}{}
    \mor{RVT}{PC}{}
    \mor{PC}{RVT}{}
    \mor{SGV}{PC}{}[+1,1]
    \enddc
    \]
    
\medskip     
    
    Here, \textit{PC} represents the Puisuex characteristic of a plane curve, \textit{RVT} represents the RVT code of a point in the Monster Tower, and \textit{SGV} represents the small growth vector of a Goursat germ.
       The arrow \{RVT\} $\longrightarrow$ \{SGV\} was given in \cite{J}, the arrow \{RVT\} $\longleftrightarrow$ \{PC\} was given in \cite{montgomerypoints},  and the arrow \{SGV\} $\longrightarrow$ \{RVT\} was given in \cite{Mo1}.

Now suppose we are given a Goursat germ whose derived vector 
is
\[ der=(\underbrace{M_1, \ M_1, \dots, M_1}_{m_1}, \ \underbrace{M_2, \ M_2, \dots, M_2}_{m_2}, \dots, \underbrace{M_{v+1}, \ M_{v+1}, \dots, M_{v+1}}_{m_{v+1}}),  \]
with $M_1<M_2<\cdots<M_v<M_{v+1}$.
Consider the set $S=\{ M_i | \ M_{i-1} \ \text{divides}\ M_i \}$.  Let $g=|S|$.  For $1\leq j\leq g$, let $N_1, N_2, \dots, N_g$ denote the elements of $S$ in decreasing order.  We always have $N_g=M_2$, since $M_1=1$. For $1\leq j\leq g$ let $M_{k_j}=N_j$.  

\begin{theorem}[\cite{Shanbrom}]
\label{maintheorem}
The corresponding Puiseux characteristic is
$ [\lambda_0; \lambda_1, \dots, \lambda_{g}] $
where
\begin{align*}
\lambda_0&= M_{v+1} \\
\lambda_j &=\sum_{i\geq k_j} m_i M_i+M_{k_j} +M_{k_j-1}
\end{align*}
for $1\leq j\leq g$.
\end{theorem}

\begin{example}
Suppose $der=(1, 1, 2, 2, 2, 2, 2, 2, 4, 6, 6, 6, 18, 24, 24)$. Note that $\lambda_0=M_{v+1}=M_{6}=24$.  We also have
$S=\{18, 4, 2\}$, and therefore $g=3$.  Then write $S=\{18, 4, 2 \}=\{N_{1}, N_2, N_3 \}=\{M_5, M_3, M_2 \}$ so that $k_1=5,\ k_2=3,$ and $k_3=2$.
Finally, we compute
\begin{align*}
\lambda_1&= \sum_{i\geq 5} m_iM_i+M_5+M_4=90 \\
\lambda_2&= \sum_{i\geq 3} m_iM_i+M_3+M_2=94 \\
\lambda_3&= \sum_{i\geq 2} m_iM_i+M_2 +M_1=103. \\
\end{align*}
The Puiseux characteristic is thus
$ [24; 90, 94, 103]. $
\end{example}

%%%%%%%%%%%%%%%%%%%%%%%%%%%%%%%%%

%\subsection{Some remarks on the isotropy method} 

%%%%%%%%%%%%%%%%%%%%%%%%%%%%%%%%%

\subsection{Spelling rules} 
%{\color{blue} Corey: is there a leap of logic here? How does one compute the Milnor number from a parametrized curve. Which formulas are you actually using that connect Milnor and other discrete invariants? I don't see the reference of Wall here.}
%The complete spelling rules for the $\R^{3}$-Semple Tower were proposed in \cite{castro4} . 
The $RVT$ code was first studied for the 
$\R^{2}$-Semple Tower and is a word in the letters $R$, $V$, and $T$ subject to a simple set of spelling rules (\cite{montgomerypoints}).  The spelling rules come from the number of critical directions 
that appear in the rank $2$ distribution that exist above each point in the planar tower.  In \cite{castro}, Montgomery, Howard and Castro began studying the $\R^{3}$-Semple Tower and extended the alphabet for the $RVT$ 
coding system to include the letters $T_{i}$ for $i = 1,2$ and $L_{j}$ for $j = 1, 2, 3$ which come from the critical planes that exist within the rank $3$ distributions at each level of the $\R^{3}$-Semple 
Tower.  The first spelling rules were obtained in \cite{castro2}, \cite{castro3}, and in \cite{castro4} we investigated the behavior of these critical planes and completed the spelling rules.  
These spelling rules for the spatial tower are given by the following result, where the ``$:$'' denotes which letters can be placed after a given 
letter.  For example, given the letter $R$ one can put either the letters $R$ or $V$ after it.  

%{\color{magenta}Above, may want to change Monster to Semple or at least mention the terminology change.  Below, in statement of theorem, should say $\mathbb R^3$ somewhere.}

\begin{theorem}{\cite{castro4}}
\label{thm:comspell}
The complete spelling rules for any $RVT$ code 
in the $\mathbb R^3$-Semple Tower
are as follows:
\

\no $(1)$ Any $RVT$ code string must begin with the letter $R$.
\

\no $(2)$ $R:$ $R$ and $V$.
\

\no $(3)$ $V$ and $T (= T_{1}):$ $R$, $V$, $T$, and $L$.
\

\no $(4)$ $L (= L_{1})$ and $L_{j}$ for $j = 2,3$: $R$, $V$, $T_{i}$ for $i = 1,2$, and $L_{j}$ for $j = 1,2,3$.
\

\no $(5)$ $T_{2}: R, V, T_{2},$ and $L_{3}$.
\end{theorem}      

The significance of this result is the role it plays in the classification problem of the points within the spatial tower.  In \cite{castro2} we used a technique called the \textit{isotropy method} which allows
us to classify points at any level of the spatial tower so long as we know how to describe the $RVT$-classes in Kumpera-Ruiz coordinates.

%{\color{magenta}Should be Ruiz not Rubin above, right?}
\begin{theorem}[\cite{castro3}]
\label{thm:orbits}
In the spatial Semple Tower the number of orbits within each of the first four levels of the tower are as follows:
\

$\bullet$ Level $1$ has $1$ orbit,
\

$\bullet$ Level $2$ has $2$ orbits,
\

$\bullet$ Level $3$ has $7$ orbits,
\

$\bullet$ Level $4$ has $34$ orbits.
\end{theorem}                                                                

%\noindent {\color{blue} Wyatt: Can you include a brief paragraph on the orbits we already counted for levels 1 to 4? }
%\
%\noindent {\color{cyan} Alex: Done.}

%%%%%%%%%%%%%%%%%%%%%%%%%%%%%%%%%

\section{Bridge 2: nonholonomic geometry and algebraic topology} 

\subsection{Semple meets Chern: nontrivality of the $\C ^{2}$-Semple Tower}

Some interesting algebraic topological questions arise when one starts to consider the complexified version of the Semple Tower. If we replace $S(0)$ by $\C^2$, the fibers over the origin cease to be a 
trivial Cartesian product as the following computation with cohomology classes show. 

We will show, for a base consisting of a neighborhood $U$ of the origin in $\C ^{2}$ that the $n$th level of the Semple Tower is \textit{not} the product manifold $U \times (\C 
P^{1})^{n} = U \times \C P^{1} \times \cdots \times \C P^{1}$.  We can show this by using the \textit{Borel-Hirzebruch formula}, found in \cite{kuroki}, in order to compute the cohomology of the 
$\C ^{2}$-Semple Tower. 
\ 

Let $\xi$ be a rank $n$ complex vector bundle over a topological space $X$, and let $\Pj (\xi)$ denote its projectivization.  Then the \textit{Borel-Hirzebruch formula} is given by      
$$H^{\ast}(\Pj (\xi); \Z) \simeq H^{\ast}(X; \Z)[x] \, / \, <x^{n} + \sum ^{n}_{i = 1} (-1)^{i} c_{i}(\pi ^{\ast} \xi)x^{n-i} )> ,$$
\no where $\pi ^{\ast} \xi$ is the pullback of $\xi$ along $\pi: \Pj (\xi) \rightarrow X$ and $c_{i}(\pi ^{\ast} \xi)$ if the $i$th Chern class of $\pi ^{\ast} \xi$.  Here, $x$ can be viewed as the first Chern class
of the canonical line bundle over $\Pj(\xi)$.  One can also replace $c_{i}(\pi ^{\ast} \xi)$ with $c_{i}(\xi)$ since the induced homomorphism $\pi ^{\ast} : H^{\ast}(X; \Z) \rightarrow H^{\ast}(\Pj (\xi) ; \Z)$ 
is injective.
We can apply the Borel-Hirzebruch formula to an $n$-level $\C P$-tower 
$$ S(m) \overset{\pi_{m}}{\longrightarrow} S(m-1) \overset{\pi _{m-1}}{\longrightarrow} \cdots \overset{\pi_{2}}{\longrightarrow} S(1) \overset{\pi_{1}}{\longrightarrow} S(0) = \{ \text{a point} \}, $$
\no with $S(i) = \Pj (\xi _{i-1})$, to get the isomorphism
$$ H^{\ast}(S(m); \Z)  \simeq \Z [x_{1}, \cdots , x_{m}] / <x^{n_{k} }_{k} + \sum^{n_{k}}_{i = 1} (-1)^{i} c_{i}(\xi _{k-1}) x^{n_{k} - i}_{k} \, | \, k = 1, \cdots , m > $$

%{\color{magenta}I don't understand the $<>$ notation in these formulas.    Is this used for set notation, like \{ \}?  It appears three times.  Also, above, how is $S(0)$= a point?  I thought $S(0)=U$.  By = do you mean contractible?  Also, note the Monster language here.} \\

%{\color{blue} From Alex: I changed Wyatt's old notation $C_m$ to be consistent with my notation for the Semple tower.} 

\begin{theorem}
For $n \geq 2$, the $n$th level of the $\C ^{2}$-Semple Tower is a nontrivial bundle and the cohomology at each level of the tower is of the form
$H^{\ast}(S(n); \Z) \simeq \Z [x_{1}, \cdots, x_{n}] / <x^{2}_{1}, x^{2}_{k} - c_{1}(\Delta_{k-1})x_{k} \, | \, k = 2, \cdots, n > $.

\begin{proof}
The first level of the Semple Tower is a trivial bundle given by $S(1)= U \times \C P^{1}$ with $U$ being a contractible open subset of the origin in $\C ^{2}$.  Our rank $2$ distribution
over $S(1)$ is $\Delta_{(p, \ell)} = d \pi ^{-1}_{(p, \ell)} (\ell)$ for $p \in U$ and $\ell \subset T_{p} \C ^{2}$.  We use the approach given in \cite{milnor} to determine the
first and second Chern classes for $\Delta_{1}$.  We note that there is a nonvanishing section $s: S(1) \rightarrow \Delta_{1}$ given by $(p, \ell) \mapsto \ell$, since $\ell$ is never zero and
hence tells us that the second Chern class of $\Delta_{1}$ will vanish.  Let $\Delta ^{0}_{1}$ be the rank $1$ subdistribution of $\Delta_{1}$ defined by 
$\Delta ^{0}_{1} (p, \ell) = \Delta _{1}(p, \ell) / span \{ \ell \}$ which will be the tangent space to $\C P^{1}$.  It is well know that $T \C P^{1}$ has nontrivial first Chern class.  This implies 
$c_{2}(\Delta_{1}) = 0$ and $c_{1}(\Delta_{1}) \neq 0$, and since $S(2) = \Pj (\Delta_{1})$ we end up with $H^{\ast}(S(2); \Z) \simeq \Z [x_{1}, x_{2}]/ <x^{2}_{1}, x^{2}_{2} - c_{1}
(\Delta_{1})x_{2} >$.
\
One can see that we can apply the same reasoning as above to show $c_{2}(\Delta_{i}) = 0$ and $c_{1}(\Delta_{i}) \neq 0$ for $i \geq 2$ and that
$H^{\ast}(S(n) ; \Z ) \simeq \Z[x_{1}, \cdots, x_{n}] / <x^{2}_{1}, x^{2}_{i} - c_{1}(\Delta_{i-1})x_{i} \, | \, i = 2, \cdots , n >$.
\end{proof}
\end{theorem}         

The main open question here is: {\em  Can one realize these cohomology classes as singularity classes within the Semple Tower?} R. Thom, who to our knowledge, was the first to propose this sort of program in algebraic topology (\cite{thom}). Whether this realization has direct applications in controllability or stabilization questions of the underlying control system remains elusive to us.

%%%%%%%%%%%%%%%%%%%%%%%%%%%%%%%%%

\subsection{Semple meets Milnor}

%\noindent {\color{blue} Corey: I shortened this session a bit. Let me know what you think. I started with your second paragraph \dots} 
%{\color{magenta} Fine.  I reworded a bit so that both paragraphs didn't start with 'In [17]'.  I think the stuff you deleted (explanation of RVT codes) could appear in Section 2.3.}

%In the case $S(0)$ has dimension 2 the RVT words are required to satisfy the following spelling rules: the letter T cannot follow the letter R.  An RVT code represents an equivalence class of Goursat germs.  The construction of a Goursat germ's RVT code was implicit in \cite{montgomerygeo}, and made explicit in \cite{Mo1}, where the letters G,S,T were used instead of R,V,T.  Beginning with a Goursat germ, one forms the Goursat flag, which possesses a canonical integrable subflag called the \emph{characteristic foliation} or \emph{Cauchy characteristic}.  The geometric relationship between the members of the two flags can be characterized as Regular, Vertical, or Tangent (or, alternatively, Generic, Singular, or Tangent) and one encodes this information into a word called the RVT code.  %See Section 1.2 of \cite{Mo1} for details.
A parallel definition of the RVT codes for Goursat germs was proposed in \cite{montgomerypoints} using 
the Semple Tower.
%a tower of manifolds called the ``Monster Tower.''  
This tower is Goursat universal: every Goursat germ occurs somewhere within the tower.  Each point in the Semple Tower is assigned an RVT code, and the code of a Goursat germ at a reference point $p$ is that of $p$ itself.  See \cite{montgomerypoints} for details, or Section 2.2 of \cite{Shanbrom} for a summary.

In \cite{montgomerypoints}, a correspondence between points in the Semple Tower and plane curve germs was made explicit.  Singular curves correspond to points whose RVT code ends with the letter V or T.  The Milnor number is a fundamental invariant of such curve singularities.  Our current work seeks to compute the Milnor number $\mu$ from a given RVT code, and we present some preliminary results below, after recalling the definition of $\mu$.

Suppose $C$ is the germ at $O$ of the singular plane curve defined by $f(x,y)=0$. Let $B_{\epsilon}$ denote the disk of radius $\epsilon$ centered at the origin in $\mathbb C^2$, with boundary sphere $S_{\epsilon}$.

\begin{definition}%[Milnor]
Let $K=f^{-1}(0) \cap S_{\epsilon}$.  Then for sufficiently small $\epsilon$, the map
\begin{align*}
\phi \colon S_{\epsilon}-K &\to S^1 \\
z &\mapsto f(z)/|f(z)|
\end{align*}
is a fibration, known as the \emph{Milnor fibration}.
\end{definition}

%\begin{definition}[Wall]
% Let $D^*_{\eta}$ denote the disk of radius $\eta$ centered at the origin in $\mathbb C$, punctured at 0. Then for sufficiently small $\epsilon$, there is an $\eta_0$ such that for $\eta<\eta_0$, the restriction of $f$
%\[f \colon B_{\epsilon} \cap f^{-1}(D^*_{\eta}) \to D^*_{\eta} \]
%is a fibration, known as the \emph{Milnor fibration}.
%\end{definition}

The fiber $F$, known as the \emph{Milnor fiber}, is a compact, connected, oriented surface with $r$ boundary components, where the curve $C$ has $r$ branches.
The first Betti number of $F$ is the \emph{Milnor number} of the singularity, denoted $\mu$.

%If $f$ has no repeated factor, then $O$ is an isolated point of the intersection $\{(x,y)\ |\ \partial f/\partial x= \partial f/ \partial y =0\}$.  This leads to the algebraic definition of the Milnor number \dots
%{\color{blue} Our starting point is an exercise (6.7.2) that appears in Wall's book (\cite{wall}). }
The two formulas below are {conjectured to} give the Milnor number for a prescribed RVT code.  The first formula (\ref{formula1}) concerns a single block of the form $R^sV^kT^u$ where the parameters $s,k$, and $u$ are arbitrary non-negative integers.  Any RVT code consists of a sequence of such blocks.  Here the superscripts denote multiplicities of letters. The second formula (\ref{formula2}) concerns RVT codes which consist of strings of the form $R^{s_j}VT^{u_j}$, where the parameters $s_j$ and $u_j$ are positive integers.  Proofs of these formulas will appear in a forthcoming paper of Howard and Shanbrom. Here $F(k)$ denotes the $k$th Fibonacci number, where $F(1)=F(2)=1$.  Also, we consider $\mu/2$ instead of the Milnor number $\mu$ for convenience, and $\begin{pmatrix} k \\ 2 \end{pmatrix}$ denotes $k$ choose 2.

%We extensively use Proposition  6.3.2 from \cite{W}, and consider $\mu/2$ instead of the Milnor number $\mu$ for convenience.  We denote the $k$th Fibonacci number by $F(k)$, where $F(1)=F(2)=1$.

 Basic building block:
\begin{equation*}\tag{*} \label{formula1}
\begin{split}
\frac{\mu}{2}(R^sV^kT^u)&=
 \begin{pmatrix} F(k+2) \\ 2 \end{pmatrix} (2+2u)
-\begin{pmatrix} F(k+1) \\ 2 \end{pmatrix}u 
+\begin{pmatrix} F(k+2)+F(k)u \\ 2 \end{pmatrix}(s-2)\\
%+ 2u \begin{pmatrix} F(k+2) \\ 2 \end{pmatrix}
&+ F(k)F(k+2)\frac{u(u+1)}{2}
+ \sum_{j=1}^{k-1}  \begin{pmatrix} F(j+2) \\ 2 \end{pmatrix}.
\end{split}
\end{equation*}
For example, $\mu(R^3V^5T^2)=1804$.

 Iterative process for single $V$'s:
\begin{equation*}\tag{**} \label{formula2}
\begin{split}
\frac{\mu}{2}(R^{s_1} V T^{u_1} R^{s_2} V T^{u_2} \cdots R^{s_n} V T^{u_n}) &= \begin{pmatrix} (u_1+2) \cdots (u_n+2) \\ 2 \end{pmatrix}s_1\\
&+\begin{pmatrix} (u_2+2) \cdots (u_n+2) \\ 2 \end{pmatrix}(s_2+u_1+1) +\dots\\
&+\begin{pmatrix} (u_j+2) \cdots (u_n+2) \\ 2 \end{pmatrix}(s_j+u_{j-1}+1)
+\dots\\
&+\begin{pmatrix} (u_n+2) \\ 2 \end{pmatrix}(s_n+u_{n-1}+1).
\end{split}
\end{equation*}
For example, $\mu(R^2VT^2 R^3VT RVT^3)=8400$.

%{\color{cyan} Maybe try to fit this equation onto just $2$ lines so it gives us a little bit more space?  Not necessary, though.}

%%%%%%%%%%%%%%%%%%%%%%%%%%%%%%%%%%%%%%%%%%%%%
%Mention all people who have directly or indirectly collaborated to the manuscript
\section*{Acknowledgements.}

The authors thank Gary Kennedy (Ohio State, Mansfield) and Susan Colley (Oberlin College) for various useful and illuminating discussions, and for keeping us busy thinking about the beautiful analogies between algebraic geometry and the local geometry of nonholonomic plane fields. Warm thanks also to M. Weissman (Yale-NUS) for first proposing to us the questions about the global topology of the Semple Tower and the puzzles related to the representation theory of jets of diffeomorphisms. Very warm thanks to David Mart\'in de Diego (ICMAT at Madrid) and Darryl Holm (Imperial College at London): the former for his hospitality in Madrid and the latter for helping us promote our work. We also would like to express here our gratitude for our teachers Richard Montgomery and Misha Zhitomirskii for pointing out the trail head that led us to this fascinating world of subriemannian geometry. % Anyone else we would like to thank?

%%%%%%%%%%%%%%%%%%%%%%%%%%%%%%%%%%%%%%%%%%%%%
%The bibliography 

\medskip
% The data information below will be filled by AIMS editorial staff
Received xxxx 20xx; revised xxxx 20xx.
\medskip

%%%%%%%%%%%%%%%%%%%%%%%%%%%%%%%%
% Appendix 

%\begin{appendices}
%  \chapter{Geometric model for the airport luggage cart}
% A quick derivation. 
%\end{appendices}

%\vspace{.25in}


\begin{thebibliography}{99}

\bibitem{arnold1} 
	\newblock  V. I. Arnold,
	\newblock  {Simple singularities of curves}, 
	\newblock  \emph{Proc. Steklov Inst. Math.}, {\bf 226} (1999), 20--28.

	
\bibitem{arnold3}
        \newblock V. I. Arnold, A.N. Varchenko, and S. Guzein-Sade,
        \newblock \emph{Singularit\'es des Applicants Diff\'erentiables},
        \newblock Editions MIR, Moscow, 1986.

\bibitem{kennedy1} 
	\newblock G. Kennedy and S. Colley,
	\newblock Triple and quadruple contact of plane curves,
	\newblock  \emph{Enumerative Algebraic Geometry} (Copenhagen, 1989), 31?59, 
Contemp. Math., 123, Amer. Math. Soc., Providence, RI, 1991.

\bibitem{shibuyayamaguchi}
       \newblock K. Shibuya and K. Yamaguchi,
       \newblock Drapeau theorem for differential systems,
       \newblock \emph{Differential Geom. Appl.}, \textbf{27} (2009), 793--808.



\bibitem{castro}
       \newblock A. Castro and W. Howard,
       \newblock A Semple-type approach to a problem of Goursat: the multi-flag case,
       \newblock \emph{C. R. Math. Acad. Sci. Paris}, \textbf{351} (2013), 921--925.
       
\bibitem{castro2}
      \newblock A. Castro, R. Montgomery, and appendix by W. Howard,
      \newblock Spatial curve singularities and the Monster/Semple tower,
      \newblock \emph{Israel J. Math.}, \textbf{192} (2012), 381--427.

\bibitem{castro3}
      \newblock A. Castro and W. Howard,
      \newblock A Monster tower approach to Goursat multi-flags,
      \newblock \emph{Differential Geom. Appl.}, \textbf{30} (2012), 405--427.

\bibitem{castro4}
      \newblock A. Castro and W. Howard,
      \newblock Spelling rules for the Monster/Semple tower,
      \newblock arXiv:1407.1824.
      
\bibitem{castro5} % results on desingularization of curves and the Enriques formula
	\newblock A. Castro, 
	\newblock \emph{Chains and Monsters}, 
	\newblock Ph.D thesis, UC, Santa Cruz, 2010.
	
\bibitem{lejeune} 
	\newblock Monique Lejeune-Jalabert,
	\newblock Chains of points in the Semple tower,
	\newblock \emph{Am. J. Math.}, \textbf{128} (2006),  1283--1311.
	
	
\bibitem{giarokumperaruiz}
        \newblock A. Giaro, A. Kumpera, and C. Ruiz,
        \newblock Sur la Lecture correcte d'un r\'esultat d' \'Elie Cartan,
        \newblock \emph{C.R. Acad. Sci. Paris S\'er. A-B}, \textbf{287} (1978), A241--A244.


	
\bibitem{hausmannrodriguez}
        \newblock J.C. Hausmann and E. Rodriguez,
        \newblock Holonomy orbits of the snake charmer algorithm,
        \newblock \emph{Geometry and Topology of Manifolds}, \textbf{76} (2007), 207--219.

\bibitem{kuroki}
        \newblock S. Kuroki and D. Suh,
        \newblock Classification of complex projective towers up to dimension 8 and cohomological rigidity,
        \newblock arXiv:1203.4403.

	
\bibitem{milnor}
        \newblock J. Milnor and J. Stasheff,
        \newblock \emph{Characteristic Classes},
        \newblock	Princeton University Press, 1974.
	
\bibitem{murray}
         \newblock R. Murray and S. Sastry,
         \newblock Nonholonomic motion planning: steering using sinusoids,
         \newblock \emph{IEEE Trans. Automat. Control}, \textbf{38} (1993), 700--716.


	
\bibitem{montgomerygeo}
         \newblock R. Montgomery and M. Zhitomirskii,
         \newblock Geometric appraoch to Goursat flags,
         \newblock \emph{Ann. Inst. H. Poincar\'e Anal. Non Lin\'eaire}, \textbf{18} (2001), 459--493.
         
\bibitem{montgomerypoints}
         \newblock R. Montgomery and M. Zhitomirksii,
         \newblock Points and curves in the Monster tower,
         \newblock \emph{Mem. Amer. Math. Soc.}, \textbf{203} (2009), x+137.


\bibitem{pelletier}
        \newblock F. Pelletier and M. Slayman,
        \newblock Configuration of an articulated arm and singularities of special multi-flags,
        \newblock \emph{SIGMA}, \textbf{10} (2014), 1--38. 



\bibitem{thom}
        \newblock R. Thom,
        \newblock Quelques propri\'et\'es globales des vari\'et\'es,
        \newblock \emph{Comment. Math. Helv.}, \textbf{28} (1954), 17--86.



\bibitem{J} 
	\newblock F. Jean,
	\newblock {The car with N trailers: characterisation of the singular configurations},
	\newblock \emph{Control, Optimisation, and Calculus of Variations}, \textbf{1} (1996), 241--266.


\bibitem{Mo1} P. Mormul,
	\newblock Geometric classes of Goursat flags and their encoding by small growth vectors,
	\newblock \emph{Central European J. Math.}, \textbf{2} (2004), 859--883.


\bibitem{Shanbrom} C. Shanbrom, 
	\newblock The Puiseux characteristic of a Goursat germ,
	\newblock \emph{J. Dynamical and Control Systems},
	\newblock\textbf{20}  (2014), 33--46.


\bibitem{W} C. Wall, 
	\newblock \emph{Singular Points of Plane Curves},  
	\newblock London Mathematical Society Student Texts, 2004.
 
\end{thebibliography}
\end{document}